\documentclass[11pt]{amsart}
\usepackage{amssymb}
\usepackage{epsfig}

\newcommand{\sevafig}[3]{\begin{figure}[h]\centerline{
 \epsfig{file=#1,width=#2,angle=#3}}
\bigskip\caption{}\end{figure}}

\newcommand{\nw}{\bigwedge\nolimits}
\newcommand{\Star}{\matho{Star}}
\newcommand{\med}{\,|\,}
\newcommand{\Conf}{\matho{Conf}\nolimits}
\def\H{\mathcal H}
\def\Imn{\matho{Im}}

\newcommand{\Arg}{\matho{Arg}}
\newcommand{\Log}{\matho{Log}}

\def\L{{L_{\infty}}}
\def\p{{\partial}}
\def\ra{{\longrightarrow}}

\newcommand{\U}{\mathcal U}
\newcommand{\F}{\mathcal F}

\newcommand{\lin}{{\mathrm{lin}}}
\newcommand{\dg}{dg~}

\renewcommand{\phi}{\varphi}

\newcommand{\ndot}{\bullet}
\newcommand{\epr}{\qed}

\newcommand{\g}{{\mathfrak g}}

\def\matho#1{\mathop{\mathrm{#1}}}

\DeclareMathSizes{11.1}{10}{8}{6}

\newcommand\nfrac[2]
{\dfrac{\raisebox{-1pt}{\fontsize{11.1}{10pt}\selectfont$#1$}}
       {\raisebox{ 1pt}{\fontsize{11.1}{10pt}\selectfont$#2$}}}

\newcommand{\Hom}{\matho{Hom}}

\newcommand{\poly}{{\mathrm{poly}}}

\newcommand{\HKR}{{\mathrm{HKR}}}

\newcommand{\C}{\mathbb C}
\newcommand{\R}{\mathbb R}
\newcommand{\Z}{\mathbb Z}

\newtheorem*{theorem}{Theorem}

\newtheorem*{lemma}{Lemma}

\theoremstyle{remark}

\theoremstyle{definition}
\newtheorem*{defin}{Definition}

\author{Boris Shoikhet}
\title[On the Kontsevich and CBH-quantizations]%
{On the Kontsevich and the Campbell--Baker--Hausdorff
deformation quantizations of a linear Poisson structure}
\date{08.03.1999}
\address{IUM, 11 Bol'shoj Vlas'evskij per.,
Moscow 121002, Russia}
\email{borya@mccme.ru}

\begin{document}
\maketitle

\sloppy

\begin{abstract}

For the Kirillov--Poisson structure on the vector space $\g^*$,
where $\g$ is a finite-dimensional Lie algebra, it is known at
least two canonical deformations quantization of this structure:
they are the M.\,Kontsevich universal formula [K], and the formula,
arising from the classical Campbell--Baker--Hausdorff formula
[Ka]. It was proved in~[Ka] that the last formula is exactly the
part of Kontsevich's formula consisting of all the admissible
graphs without (oriented) cycles between the vertices of the
first type. It follows from the CBH-theorem that this part of
Kontsevich's formula defines an \emph{associative} product
(in the case of a linear Poisson structure).

The aim of these notes is to prove the last result directly, using
the methods analogous to [K] instead of the CBH-formula. We
construct an $L_\infty$-morphism $\U_\lin\colon[T^\ndot_\poly]_\lin\to
D^\ndot_\poly$ from the \dg Lie algebra of polyvector fields with
\emph{linear coefficients} to the \dg Lie algebra of
polydifferential operators, which is \emph{not} equal to the
restriction of the Formality \hbox{$L_\infty$-morphism}
$\U\colon T^\ndot_\poly\to D^\ndot_\poly$ [K] to the subalgebra
$[T^\ndot_\poly]_\lin$. For a bivector field~$\alpha$ with linear
coefficients such that $[\alpha,\alpha]=0$ the corresponding
solution $\U_\lin(\alpha)$ of the Maurer--Cartan equation in
$D^\ndot_\poly$ defines exactly the CBH-quantization,in the case of the harmonic angle 
map [K], Sect.2.We prove the associativity of the restricted Kontsevich 
formula (in the linear case) also for any angle map [K], Sect.6.2.
\end{abstract}

\section{$L_\infty$-morphisms, the Maurer--Cartan equation, and
$*$-products}

Let $\F\colon T^\ndot_\poly\to D^\ndot_\poly$ be an $L_\infty$-morphism
from the \dg Lie algebra of polyvector fields on $\R^d$ to the
\dg Lie algebra of polydifferential operators on $\R^d$, and let
\begin{gather*}
\F_1\colon T^\ndot_\poly\to D^\ndot_\poly\\
\F_2\colon\wedge^2T^\ndot_\poly\to D^\ndot_\poly [-1]\\
\F_3\colon\wedge^3T^\ndot_\poly\to D^\ndot_\poly [-2]\\
\hbox to 4cm{\dotfill}
\end{gather*}
be its Taylor components.

Than any solution $\alpha\in T^1_\poly$ of the Maurer--Cartan
equation (i.e. $\alpha$ is a bivector field such that
$[\alpha,\alpha]=0$) defines a solution $\F(\alpha)\in D^1_\poly$
of the Maurer--Cartan equation in~%
$D^\ndot_\poly$ $(\F(\alpha)\in\Hom_\C(C^\infty(\R^d)^{\otimes
2}\to C^\infty(\R^d)))$ as follows:
\begin{equation}
\F(\alpha)=\F_1(\alpha)+\nfrac12 \F_2(\alpha,\alpha)+\nfrac16
\F_3(\alpha,\alpha,\alpha)+\ldots+\nfrac 1{n!} \F_n(\alpha,\dots,\alpha)+\dots
\end{equation}

One can prove that the bidifferential operator $\F(\alpha)$ satisfy
the Maurer--Cartan equation
\begin{equation}
d\F(\alpha)+\nfrac12 [\F(\alpha),\F(\alpha)]=0
\end{equation}
where $d$ is the Hochschild differential and $[,]$ is the
Gerstenhaber bracket. It follows directly from the definitions
that $(2)$ is equivalent to the statement that the formula
\begin{equation}
f*g=f\cdot g+\F(\alpha)(f\otimes g)
\end{equation}
defines an \emph{associative} product on the vector space
$C^\infty(\R^d)$.

\section{Formality $L_\infty$-morphism $\U\colon T^\ndot_\poly\to
D^\ndot_\poly$ \rm[K]}

\subsection {Admissible Graphs and Weights}\label{s131}

\begin{defin} Admissible graph is an oriented graph with
labels such that

1) the set of vertices $V_\Gamma$ is
$\{1,\dots,n\}\coprod\{\overline1,\dots,\overline m\}$ where
$n,m\in\Z_{\ge0}$; vertices from the set $\{1,\dots,n\}$ are
called vertices of the first type, vertices from
$\{\overline1,\dots,\overline m\}$ are called vertices of the second type,

2) every edge $(v_1,v_2)\in E_\Gamma$ starts at a vertex of
first type, $v_1\in\{1,\dots,n\}$;

3) there are no loops, i.e. no edges of the type $(v,v)$;

4) for every vertex $k\in\{1,\dots,n\}$   of the first type, the
set of edges
$$\Star(k):=\{\,(v_1,v_2)\in E_\Gamma \med v_1=k\,\}.$$
starting from $s$ is labeled by symbols
$(e_k^1,\dots,e_k^{\Star(k)})$.
\end{defin}
For any admissible graph $\Gamma$, we define weight
$W_\Gamma\in\C$ by formula
$$
W_\Gamma=\prod_{k=1}^n\nfrac1{(\#\Star(k))!}\cdot
\nfrac1{(2\pi)^{2n+m-2}}\int_{C_{n,m}^+}
\bigwedge_{e\in E_\Gamma}d\varphi_e.
$$

Let us explain written here. Let
$$
\Conf_{n,m}=\{\,(p_1,\dots,p_n;q_1,\dots,q_m)\med
p_i\in \H,\ q_j\in\R,\ p_{i_1}\ne p_{i_2}\ \text{for $i_1\ne i_2$ and
$q_{j_1}\ne q_{j_2}$ for $j_1\ne j_2$}\,\}.
$$

Here $\H=\{\,z\in\C\med\Imn z>0\,\}$. Let $G$ be a group of affine
transformations
$G=\{\,z\mapsto az+b,\ a,b\in\R,\ a>0\,\}$.

Then
$\Conf_{n,m}^+=\{\,p_1,\dots,p_n;q_1,\dots,q_m)\in\Conf_{n,m}
\med q_1<q_2<\ldots<q_m\,\}$ is invariant under the action of $G$, and
we define $C_{n,m}=\Conf_{n,m}/G$,\ \ $C_{n,m}^+=\Conf_{n,m}^+/G$.

Every edge $e\in E_\Gamma$  defines a map from $\Conf_{n,m}$ to
$\Conf_{2,0}$ (if two end-points of $e$ are of the first type)
and to $\Conf_{1,1}$ otherwise. For $p,q\in \H\bigsqcup\R$\ \
$(p\ne q)$ we define function
$$
\Phi(p,q)=\Arg\left(\nfrac{(q-p)}{(q-\overline p)}\right)=
\nfrac1{2i}\Log\left(\nfrac{(q-p)(\overline q-p)}
{(q-\overline p)(\overline q-\overline p)}\right)
$$
and $1$-form  $d\Phi$.

This function is $G$-invariant, and this construction defines a
$1$-form $d\Phi_e$ for any $e\in E_\Gamma$, which is the
pull-back of $d\Phi$.

\begin{lemma} Integral in the definition of $W_{\Gamma}$ is
absolutely convergent for any $\Gamma$.
\end{lemma}

The proof is done in
Section 5 of~[K].\epr

\subsection {Formality Morphism}\label{s132}

\def\G{\Gamma}
\def\ga{\gamma}
\def\U{{\mathcal U}}

        For any admissible graph $\G$ with $n$ vertices of the first type,
         $m$ vertices of the second type, and $2n+m-2+l$ edges where $l\in \Z$,
         we define a
        linear map  $\U_\G:\otimes
        ^n T_{poly}(\R^d)\ra D_{poly}(\R^d)[1+l-n]$.
         This map has only one
          non-zero graded component
          $(\U_\G)_{(k_1,\dots,k_n)}$ where $k_i=\# Star(i)-1,\,\,
          i=1,\dots,n$. If $l=0$ then from $\U_\G$
           after anti-symmetrization
           we obtain a pre-$\L$-morphism.

            Let $\ga_1,\dots,\ga_n$ be polyvector fields on $\R^d$ of degrees
            $(k_1+1),\dots,(k_n+1)$, and $f_1,\dots, f_m$
              be functions on $\R^d$. We are going to write a
               formula
             for  function $\Phi$ on $\R^n$:
             $$\Phi:=\left(
             { \U}_\G(\ga_1\otimes\dots\otimes
             \ga_n)\right)(f_1\otimes
             \dots \otimes f_m)\,\,\,.$$

               The formula for $\Phi$ is the sum over all
                configurations of indices running from $1$ to $d$,
                labeled by $E_\G$:
                $$\Phi=\sum_{I:E_\G\ra \{1,\dots,d\}}\Phi_I\,\,,$$
              where $\Phi_I$ is the product over all $n+m$ vertices of
              $\G$ of certain partial derivatives of functions
               $g_j$ and of coefficients
               of $\ga_i$.

               Namely, with each vertex $i,\,\,1\le i\le n$ of the
               first type
                we associate function
               $\psi_i$ on $\R^d$
                which is a coefficient of the polyvector field
                $\ga_i$:
                $$\psi_i=\langle \ga_i, dx^{I(e^1_i)}\otimes\dots
                \otimes dx^{I(e^{k_i+1}_i)}\rangle\,\,\,.$$
                Here we use the identification of polyvector
                 fields with skew-symmetric tensor fields as
                 $$\xi_1\wedge\dots\wedge\xi_{k+1}\ra
                 \sum_{\sigma\in S_{k+1}} sgn(\sigma)\,
                  \xi_{\sigma_1}\otimes\dots
                  \otimes \xi_{\sigma_{k+1}
                  }\in \G(\R^d,T^{\otimes(k+1)})\,\,\,.$$
                 For each vertex $\overline j$ of second type the associated
                function $\psi_{\overline j}$ is defined as $f_j$.

             Now, at each vertex of graph $\G$ we put a function
             on $\R^d$ (i.e. $\psi_i$
             or $\psi_{\overline j}$). Also, on
             edges of graph $\G$ there are indices
              $I(e)$
              which label coordinates in $\R^d$. In the next step
               we put into each vertex
              $v$
              instead of function $\psi_v$ its partial derivative
               $$\left(\prod_{e\in E_\G,\,e=(*,v)}
               \p_{I(e)}\right) \psi_v,$$
               and then take the product over all vertices $v$ of $\G$.
               The result is by definition the summand~$\Phi_I$.

              Construction of the function $\Phi$ from the graph
                $\G$, polyvector fields $\ga_i$ and functions $f_j$,
                 is  invariant under the action of the group
                 of affine transformations of $\R^d$ because
                 we contract upper and lower indices.

We define an  $L_\infty$-morphism
$\U\colon T^\ndot_\poly(\R^d)\to D^\ndot_\poly(\R^d)$
by the formula for its $n$-th derivative $\U_n$,\ \ $n\ge1$,
considered as a skew-symmetric polylinear map 
\begin{gather*}
\U_n\colon \otimes^n T^\ndot_\poly(\R^d)
\to D^\ndot_\poly(\R^d)[1-n]:\\
\U_n=\sum_{m\ge0}\sum_{\Gamma\in
G_{n,m}}W_\Gamma\times \U_\Gamma
\end{gather*}

Here $G_{n,m}$ denotes the set of all admissible graphs with $n$
vertices of the first type, $m$~vertices of the second type and
$2n+m-2$ edges, $n\ge 1$,$m\ge0$.

\begin{theorem}[{[K]}, Sect.~6.4] $\U$ is  $L_\infty$-morphism
and  also  a   quasi-isomorphism.    $L_\infty$-morphism   $\U$   is
equivariant under
affine transformations.
\end{theorem}

The fact that $\U$ is quasi-isomorphism follows directly from the
fact that $\U_1=\varphi_{\HKR}$  and Hochschild--Kostant--Rosenberg
Theorem.

Formula (3) applied to the $L_\infty$-morphism $\U$ defines the
Kontsevich universal formula for the deformation quantization.

\section{A sketch of the proof of Theorem 2 [K]}

The condition that $\U$ is an $L_\infty$-morphism can be written
as follows:
\begin{equation}
\begin{gathered}
f_1\cdot\left(\U_n(\gamma_1\wedge\dots\wedge\gamma_n)\right)(f_2\otimes\ldots
\otimes
f_m)\pm\left(\U_n(\gamma_1\wedge\dots\wedge\gamma_n)\right)(f_1\otimes\ldots
\otimes
f_{m-1})\cdot
f_m+\\
+\sum^{m-1}_{i=1}\pm\left(\U_n(\gamma_1\wedge\dots\wedge\gamma_n)\right)
\left(f_1\otimes\ldots\otimes(f_i\cdot
f_{i+1})\otimes\ldots\otimes f_m\right)+\\
+\sum_{i\ne
j}\pm\left(\U_{n-1}([\gamma_i,\gamma_j]\wedge\dots\wedge\widehat\gamma_i\wedge\dots
\wedge\widehat\gamma_j\wedge\dots\wedge\gamma_n)\right)(f_1\otimes\ldots\otimes
f_m)+\\
+\nfrac12\sum_{k,l\ge 1,
k+l=n}\nfrac1{k!l!}\sum_{\sigma\in\sum_n}\pm[\U_k(\gamma_{\sigma_1}\wedge\dots
\wedge\gamma_{\sigma_k}),\U_l(\gamma_{\sigma_{k+1}}\wedge\dots\wedge
\gamma_{\sigma_n})](f_1\otimes\ldots\otimes
f_m)=0.
\end{gathered}
\end{equation}

It is clear that one can write the r.h.s.\  of (4) as a linear
combination
\begin{equation}
\sum_\Gamma c_\Gamma\cdot
\U_\Gamma(\gamma_1\otimes\ldots\otimes\gamma_n)(f_1\otimes\ldots\otimes
f_m)
\end{equation}
of expressions $\U_\Gamma$ for admissible graphs $\Gamma$ with $n$
vertices of the first type, $m$ vertices of the second type, and
$2n+m-3$ edges, where $n,m\ge 0$, $2n+m-3\ge 0$.

The coefficients $c_\Gamma$ are equal to the quadratic-linear
combinations of the weights $W_{\Gamma'}$.

We want to check that the coefficient $c_\Gamma$ vanishes for
each graph $\Gamma$.

The idea is to identity the coefficient $c_\Gamma$ with the
integral over the boundary $\partial\overline C_{n,m}$ of the
closed differential form constructed from $\Gamma$ as in Sect.~2.
The Stokes formula gives the vanishing:
\begin{equation}
\int_{\partial\overline C_{n,m}}\nw_{e\in E_\Gamma}d\Phi_e=\int_{\overline
C_{n,m}}d\left(\nw_{e\in E_\Gamma}d\Phi_e\right)=0.
\end{equation}

The boundary strata of codimension 1 are of the following two
types (see [K], Sect.~5):

\begin{itemize}
\item[ (S1):] points from subset $S\subset\{1,\dots,n\}$, $\#S\ge 2$ of
the first type move close to each other; the corresponding
boundary stratum is equal to $\partial_S\overline C_{n,m}=\overline
C_{\#S}\times\overline C_{n-\#S+1,m}$
\item[(S2):] points from subset $S\subset\{1,\dots,n\}$ of the first
type and points from the subset $S'\subset\{1,\dots,m\}$ of the second
type, such that $2\#S+\#S'\ge 2$, $\#S+\#S'\le n+m-1$, move close to
each other and to~$\R$; the boundary stratum is equal to
$$
\partial_{S,S'}\overline C_{n,m}=\overline C_{\#S,\#S'}\times\overline 
C_{n-\#S,m-\#S'+1}.
$$
\end{itemize}

One have:
\begin{equation}
0=\int_{\partial\overline C_{n,m}}\nw_{e\in
E_\Gamma}d\Phi_e=\sum_S\int_{\partial_S\overline C_{n,m}}\nw_{e\in
E_\Gamma}d\Phi_e+\sum_{S,S'}\int_{\partial_{S,S'}\overline
C_{n,m}}\nw_{e\in E_\Gamma}d\Phi_e.
\end{equation}

The idea is to identity the summands of the last sum with
summands of the r.h.s.\ of~(4).

\emph{Case} S1: the integral (7) vanishes except the case
$\#S=2$ and the two points are connected by an edge $\vec e$ (see
[K], Sect.~6.6). The case $\#S=2$ corresponds to the summands of
(4) with the bracket of polyvector fields. The integral
in the r.h.s.\  of (7) is equal, up to $2\pi$, to
the integral corresponded to the graph $\Gamma_1$ obtained from the
graph $\Gamma$ by the contraction of the edge $\vec e.$ Let us
note, that the graph $\Gamma_1$ has $(n-1)$ vertices of the first
type, $m$ vertices of the second type, and $2n+m-4=2(n-1)+m-2$
edges.

\sevafig{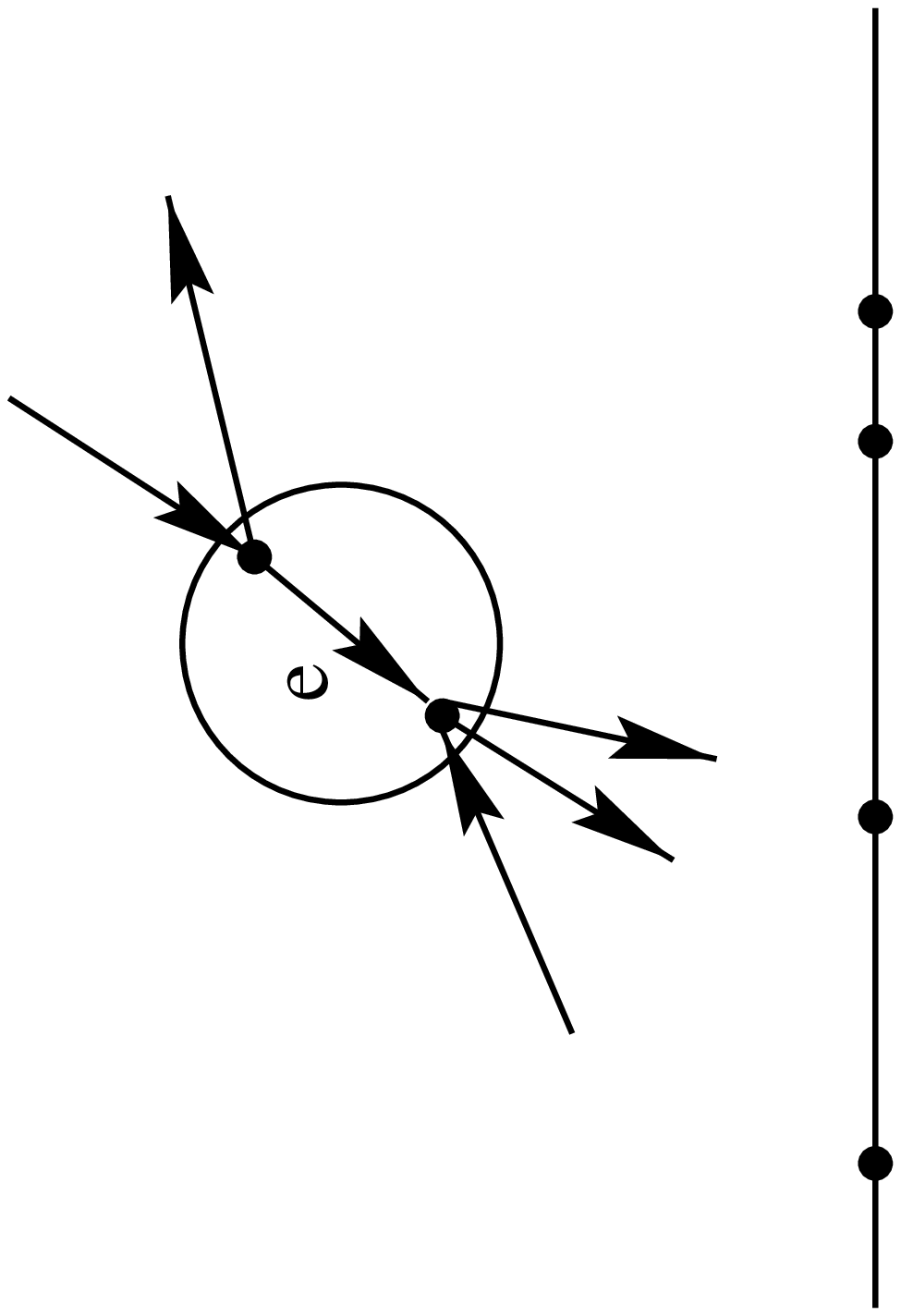}{40mm}{270}

\emph{Case} 2: One can show ([K], Sect.6.4.2.2) that the integral in (7)
vanishes except the case when there does not exist any ``external''
edge starting in the points of the subset $S\sqcup S'$. The
typical situation is shown on Fig.~2.

\sevafig{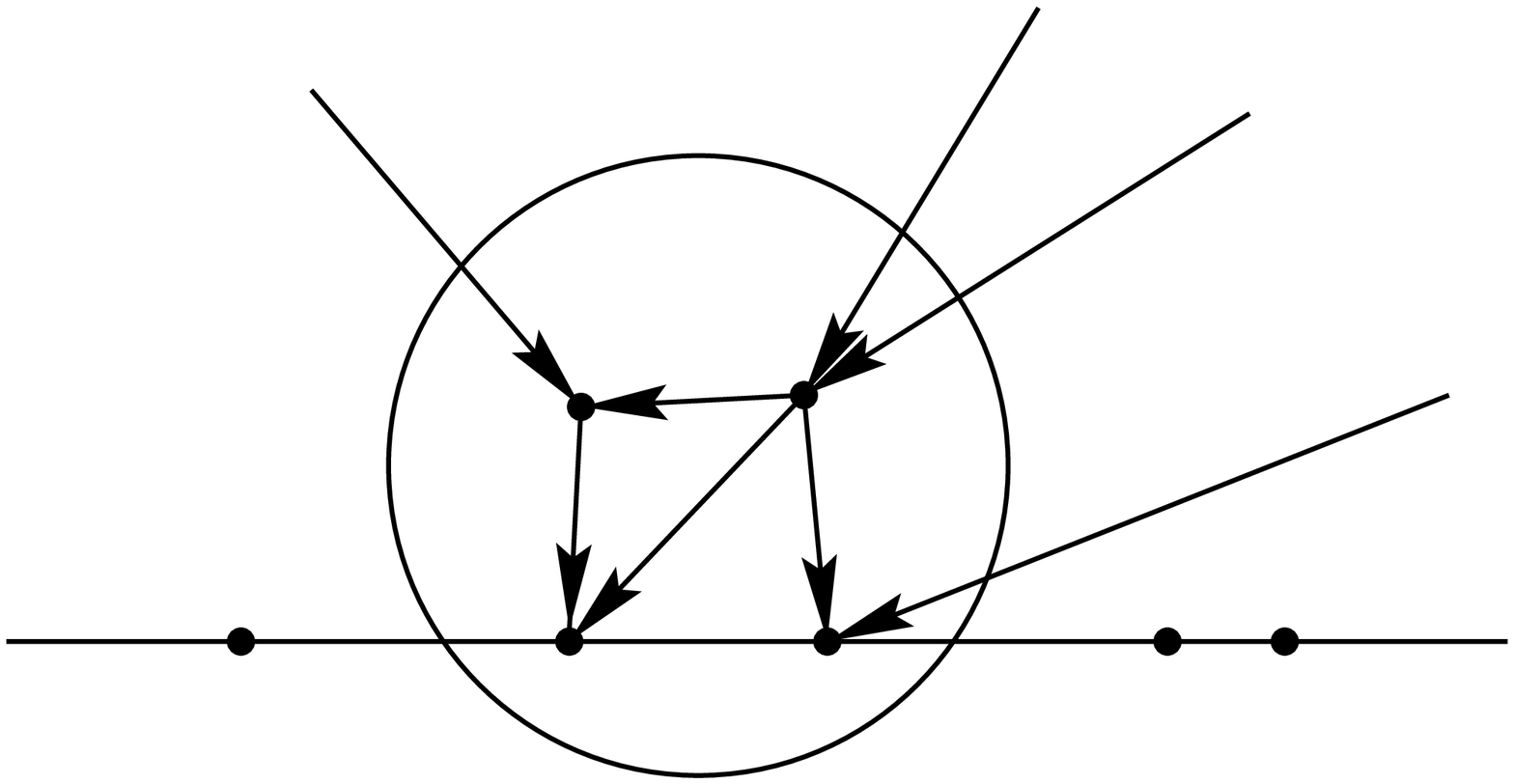}{60mm}{0}

This case is corresponded to the Gerstenhaber bracket of
polydifferential operators in the r.h.s.\  of (4). The integral is
equal to the product of the two weights $W_{\Gamma_1}\times
W_{\Gamma_2}$.

\emph{How to calculate the coefficient $c_\Gamma$}:

 Let $\Gamma$ be an admissible graph with $n$ vertices of the
 first type, $m$ vertices of the second type, and $2n+m-3$ edges.
 We consider the following two types of representations of the
 graph $\Gamma$:
 \begin{itemize}
 \item[(R1):] it is a representation of the form
 $\Gamma=\Gamma'\sqcup\vec e$, where the edge $\vec e$ connects
 two vertices of $\Gamma'$ of the first type (see Fig. 1)
 \item[(R2):] it is a representation of the form
 $\Gamma=\Gamma_1\cup\Gamma_2$ where: 1) both graphs $\Gamma_2$
 and $\overline\Gamma_1=\Gamma/\Gamma_2$ (the contraction of
 $\Gamma_2$ to a vertex of the second type) have $n_i$ vertices
 of the first type, $m_i$ vertices of the second type, and
 $2n_i+m_i-2$ edges $(i=1,2)$;\newline
 2) there does not exists any edge starting in the new vertex
 $=[\Gamma_2]$ of the graph~$\overline\Gamma_1$. (See Fig.~2).
 \end{itemize}

 Any representation of the types (R1), (R2) of the graph $\Gamma$
 has a contribution in the coefficient $c_\Gamma$, and $c_\Gamma$
 is the sum over all the possible representations. According to
 the Stokes formula, the sum of all these contributions is equal
 to 0. On the other hand, the contributions of the
 representations are in 1--1 correspondence with summands in the
 r.h.s.\  of (4).

 \section{We want to prove the $L_\infty$-Formality Conjecture for
 $\R^\infty$, or how the \dg Lie algebra $[T^\ndot_\poly]_\lin$
 appears.}

 The difficulty in the problem of the extending of the result of Section~2,~3 for
 the space $\R^\infty$ (in any sense) is the divergence of the
 polydifferential operators corresponded to the graphs with
 oriented cycles (between vertices of the first type), as is
 shown on Fig.~3.

\sevafig{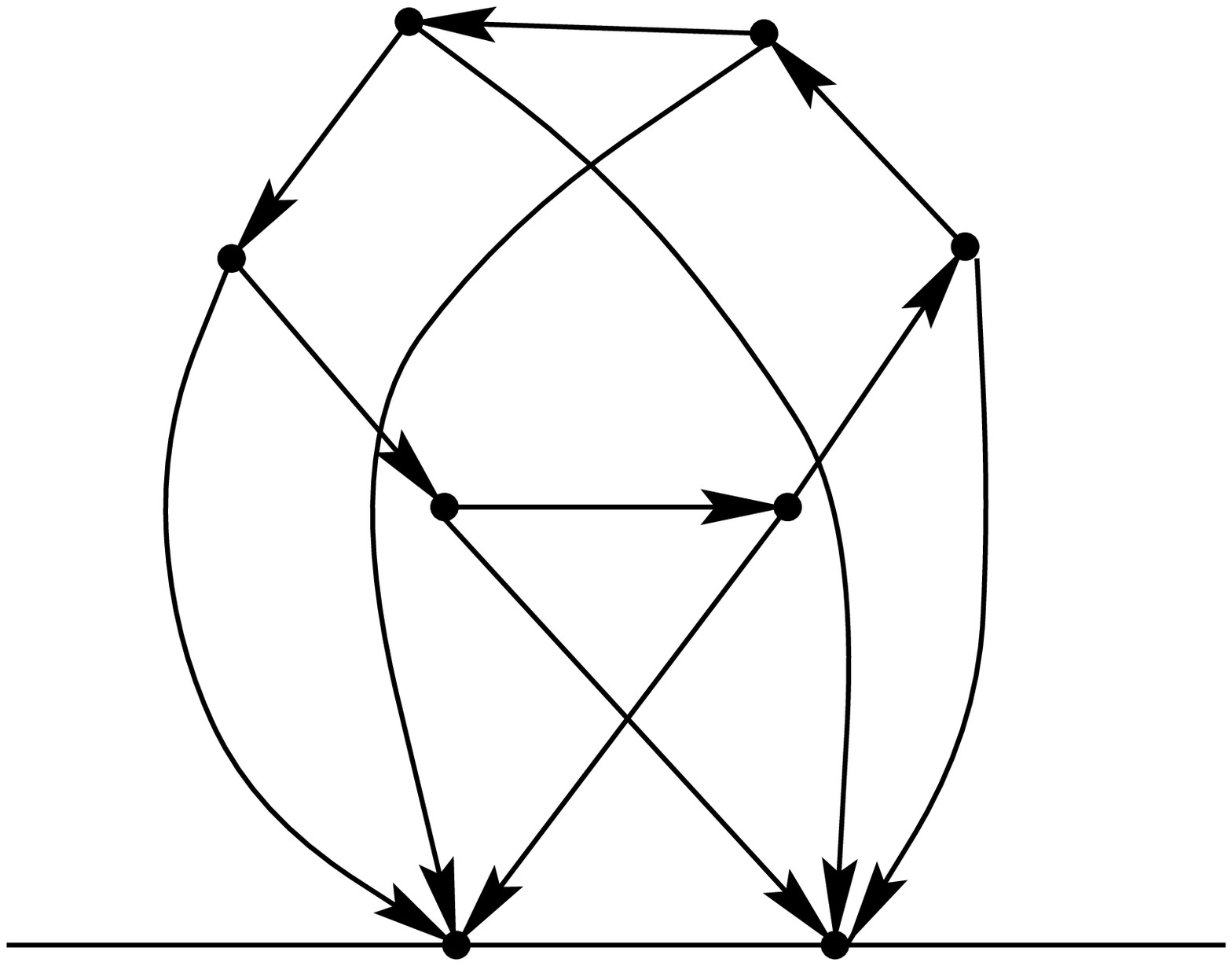}{50mm}{0}

 We want to define a new class of ``restricted'' admissible graphs
 for the definition of the $L_\infty$-morphism $\U$ in the Section
 2 such that:
 \begin{itemize}
 \item[ (i)] restricted admissible graphs do not contain any
 oriented cycles;
 \item[(ii)] the class of restricted admissible graphs is
 compatible with the two operations (R1) and (R2) (see Sect.~3),
 in the sense explained below.
 \end{itemize}

 We claim that such a class of restricted admissible graphs does
 not exist.

 \subsection{We try to exclude all the graphs with oriented
 cycles.}

 Let us suppose that the restricted class of admissible graphs
 contains graphs with non-oriented cycles between vertices of the
 first type, for example, a graph with a cycle such that all its
 edges have right orientation except the one unique edge $\vec e$
 (the general case is the same). Then the representation of the
 type (R1) $\Gamma=\Gamma'\sqcup\vec e$ does not appear in the
 r.h.s.\  of (4) because the graph $\Gamma'$ is not restricted
 admissible and is not appeared in the definition of the
 $L_\infty$-morphism $\U$. On the other hand, this representation
 appears in formula (7). Consequently, the summand in (4) and in
 (7) are not in 1--1 correspondence, and we may not use the
 arguments of the Stokes formula. Therefore, restricted
 admissible graphs may not have any (non-oriented) cycle between
 the vertices of the first type.

 \subsection{Restricted admissible graphs may not have any cycle
 (formed by vertices both of the first and the second types).}

 Let us suppose that there exists a restricted admissible graph
 containing any cycle, like is shown on the Fig. 4.
\sevafig{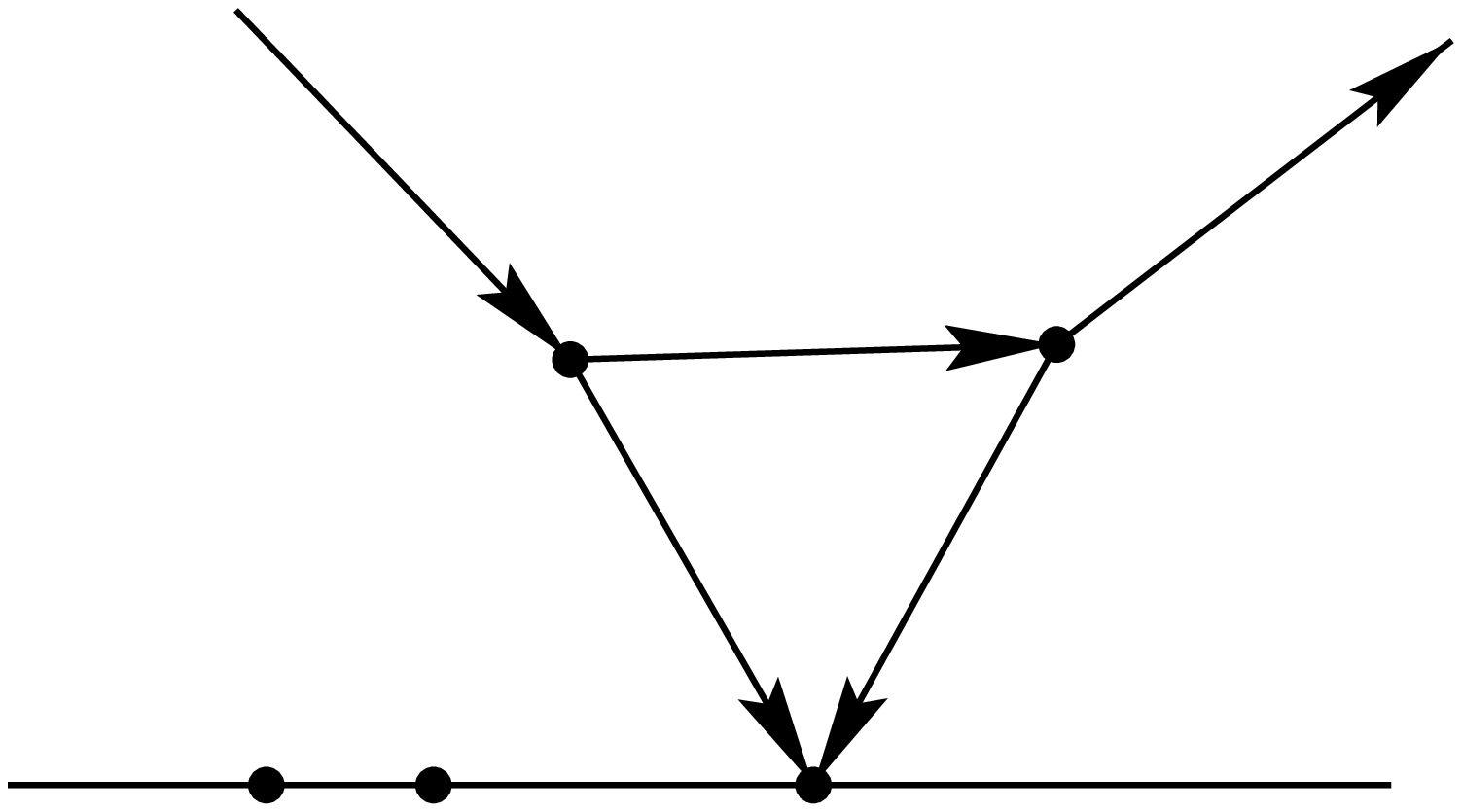}{50mm}{0}

 Then the Gerstenhaber bracket generates a graph with
 (non-oriented) cycle between vertices of the first type, as is
 shown in Fig. 5.

\sevafig{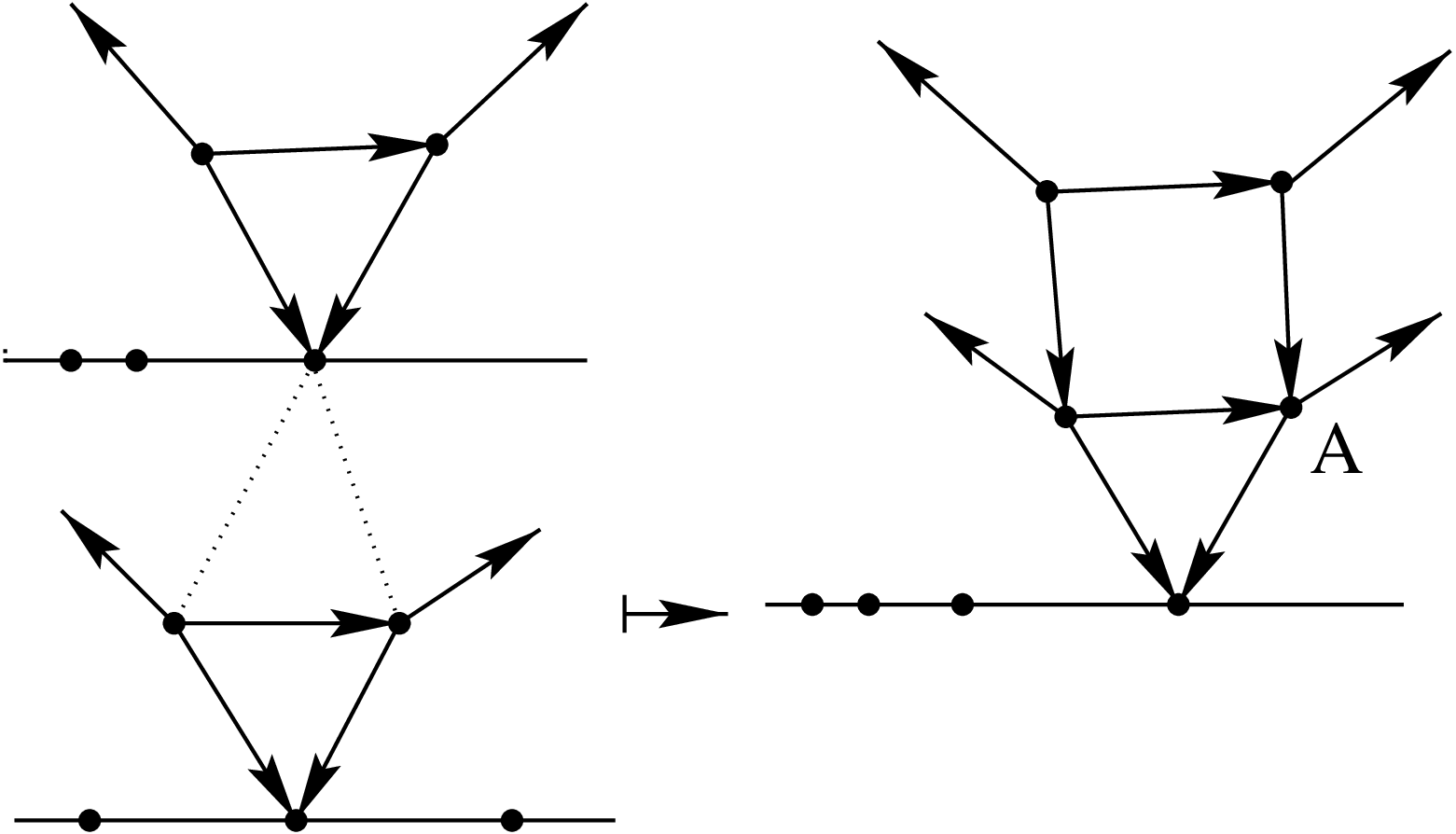}{100mm}{0}

 The graph shown at the right side of Figure 5 should be
 restricted admissible, in a contradiction with Sect.~4.l. Then,
 the restricted admissible graphs may not have any cycle, and it is
 an easy exercise to prove that the set of these graphs is empty.

 \subsection{Differential graded Lie algebra
 $[T^\ndot_\poly]_\lin$ of polyvector fields with linear
 coefficients.}

 The situation described in Sect.~4.2 will not appear for the \dg
 Lie algebra $[T^\ndot_\poly]_\lin$. Indeed, let us suppose that
 we consider only \emph{linear} polyvector fields. Then,
 there exists not more than 1 edge ending at each vertex of the
 first type. The situation of Sect.~4.2 will not appear because
 there exist 2 edges ending at the vertex $A$ on the right-hand
 side of Fig. 5, and the right-hand graph defines zero
 polydifferential operator.

 Let us summarize. Let $G^r_{n,m}$ be the set of admissible
 graphs (see Definition 2.1) with $n$ vertices of the first type,
 $m$ vertices of the second type, $2n+m-2$ edges, and which do
 not contain any (non-oriented) cycle between the vertices of the
 first type. The map
 $\U^\lin_n\colon\otimes^n[T^\ndot_\poly]_\lin\to
 D^\ndot_\poly[1-n]$ is defined as follows:
 \begin{equation}
 \U_n=\sum_{m\ge 0}\sum_{\Gamma\in G^r_{n,m}}W_\Gamma\times \U_\Gamma
 \end{equation}
 where the weight $W_\Gamma$ and the polydifferential operator
 $\U_\Gamma$ are defined as in Sect.~2. Then formula (8) defines
 the components of the $L_\infty$-morphism
 $$
 \U_\lin\colon[T^\ndot_\poly]_\lin\to D^\ndot_\poly
 $$
 which defines, by formulas (1), (3), a deformation quantization
 of the Kirillov--Poisson structure on $\g^*$ (both in
 finite-dimensional and infinite-dimensional cases). This
 deformation quantization is exactly the ``restricted''
 Kontsevich's universal formula, i.e. the Kontsevich's formula
 without graphs with any (=oriented in the linear case) cycles
 between vertices of the first type. According to the 
theorem of V.\,Kathotia [Ka], it is
 exactly the CBH-quantization.


\begin{thebibliography}{9}
 \bibitem[K]{K} M.~Kontsevich, \emph{Deformation quantization of
 Poisson manifolds}, I, preprint math.QA/9709040.
 \bibitem[Ka]{Ka} V.~Kathotia, \emph{Kontsevich's universal formula for
 deformation quantization and the Campbell--Baker--Hausdorff
 formula}, preprint math.QA/9811174 v2
 \end{thebibliography}
 \end{document}